\documentclass[11pt]{article}
\usepackage{amsmath, amssymb}
\usepackage[mathscr]{euscript} 

\usepackage[parfill]{parskip} 

\usepackage{rawfonts}
\usepackage{pictexwd}

\usepackage[latin1]{inputenc}

\newcommand{\proof}{\vspace{0pt}{\bfseries Proof.\ }} 
\newcommand{\QED}{{\unskip\nobreak\hfil\penalty50\hskip1em\hbox{}\nobreak  
   \hfil \ensuremath{\Box}\parfillskip=0pt \par}}

\newtheorem{Theorem}{Theorem}[section]
\newtheorem{Lemma}[Theorem]{Lemma}
\newtheorem{Corollary}[Theorem]{Corollary}
\newtheorem{Proposition}[Theorem]{Proposition}
\newtheorem{Remark}[Theorem]{Remark}
\newtheorem{Example}[Theorem]{Example}

\newtheorem{Definition}[Theorem]{Definition}

\usepackage[spanish]{babel}
\usepackage{endnotes}

\title{\vspace{-2.7em}  %Distancia entre borde superior de la hoja y el titulo
\Large\bfseries On well-covered, vertex decomposable and Cohen-Macaulay graphs.}
\author{Iván D. Castrill\'on\footnote{Partially supported by CONACYT and ABACUS-CINVESTAV}, Roberto Cruz and Enrique Reyes\footnote{Partially supported by SNI and ABACUS-CINVESTAV}}

\begin{document}

\maketitle
\thispagestyle{empty}

\begin{abstract}
\noindent
Let $G=(V,E)$ be a graph. If $G$ is a König graph or $G$ is a graph without 3-cycles and 5-cycle, we prove that the following conditions are equivalent: $\Delta_{G}$ is pure shellable, $R/I_{\Delta}$ is Cohen-Macaulay, $G$ is unmixed vertex decomposable graph and $G$ is well-covered with a perfect matching of König type $e_{1},...,e_{g}$ without square with two $e_i$'s. We characterize well-covered graphs without 3-cycles, 5-cycles and 7-cycles. Also, we study when graphs without 3-cycles and 5-cycles are vertex decomposable or shellable. Furthermore,  we give some properties and relations between cri\-ti\-cal, extendables and shedding vertices. Finally, we characterize unicyclic graphs with each one of the following properties: unmixed, vertex decomposable, shellable and Cohen-Macaulay.  
\end{abstract}


\section{Introduction}

Let $G$ be a simple graph (without loops and multiply edges) whose vertex set is $V(G)=\{x_{1},...,x_{n}\}$ and edge set $E(G)$. Let $R=k[x_{1},...,x_{n}]$ be a polynomial ring over a field $k$, the $edge$ $ideal$ of $G$, denoted by $I(G)$, is the ideal of $R$ generated by all monomials $x_{i}x_{j}$ such that $\{x_{i},x_{j}\}\in E(G)$. $G$ is a \textit{Cohen-Macaulay graph} if $R/I(G)$ is a Cohen-Macaulay ring (See \cite{BH}, \cite{Villa}). A subset $F$ of $V(G)$ is a \textit{stable set} or \textit{independent set} if $e\nsubseteq F$ for each $e\in E(G)$. The cardinality of the maximum stable set is denoted by $\beta (G)$. $G$ is called \textit{well-covered} if every maximal stable set has the same cardinality. On the other hand, a subset $D$ of $V(G)$ is a vertex cover of $G$ if $D\cap e \neq \emptyset$ for every $e\in E(G)$. The number of vertices in a minimum vertex cover of $G$ is called the \textit{covering} \textit{number} of $G$ and it is denoted by $\tau (G)$. This number coincide with ht$(I(G))$, the \textit{height} of $I(G)$. If the minimal vertex covers have the same cardinality then $G$ is called \textit{unmixed} graph. Notice that, $D$ is a vertex cover if and only if $V(G)\setminus D$ is an stable set. Hence, $\tau(G)=n-\beta(G)$ and $G$ is well-covered if and only if $G$ is unmixed. The \textit{Stanley-Reisner} \textit{complex} of $I(G)$, denoted by $\Delta_{G}$, is the simplicial complex whose faces are the stable sets of $G$. Recall that a simplicial complex $\Delta$ is called \textit{pure} if every facets has the same number of elements. Thus, $\Delta_{G}$ is pure if and only if $G$ is well-covered. The well-covered property has been studied for some families of graphs: graphs with girth at least 5 (in \cite{Finbow}), graphs without 4 or 5-cycles (in \cite{Finbow2}), simplicial, chordal and circular graphs (in \cite{Prisner}), block-cactus graphs (in \cite{cactus}). In this paper (Theorem \ref{24}) we characterize the well-covered graph without 3-cycles, 5-cycles and 7-cycles. 

In section 2 we give some properties and relation between critical, shedding and extendable vertices. We use these result for to study the following properties: well-covered and vertex decomposable of $G$; shellability of $\Delta_{G}$ and Cohen-Macaulay property of $R/I(G)$. These properties have been studied in (\cite{BH}, \cite{CC}, \cite{Villa5}, \cite{HH}, \cite{ON}, \cite{MM}, \cite{Stanley}, \cite{VT}, \cite{Villa}, \cite{RW}). In general, we have that (see \cite{Stanley}, \cite{Villa})

vertex-decomposable $\Rightarrow$ shellable $\Rightarrow$ sequentially Cohen-Macaulay

It is know that if each chordless cycle of $G$ has length 3 or 5 then, $G$ is vertex decomposable. (see \cite{RW}). Now, in section 3 we characterize vertex decomposable graphs without 3-cycles and 5-cycles. This result generalize the criterion given in \cite{VT} for shellable graphs. Furthermore, we characterize shellable graphs with girth at least 11.

On the other hand, we have following implications:
 
Unmixed vertex decomposable $\Rightarrow$ pure shellable $\Rightarrow$  Cohen-Macaulay $\Rightarrow$ well-covered
  
The equivalence between  Cohen-Macaulay property and unmixed vertex decomposable property has been studied for some families of graphs: bipartite graphs (in \cite{Villa5} and \cite{HH}); very well covered graphs (in \cite{half} and \cite{MM}); graphs with girth at least $5$, block-cactus and graphs without 4-cycles and 5-cycles (in \cite{ON}). In section 5 we prove the equivalences of these properties for König graphs and graphs without 3-cycles and 5-cycles. Furthermore, we prove theses properties are equivalent to the following condition: $G$ is unmixed graph with a perfect matching $e_{1}=\{x_{1},y_{1}\},...,e_{g}=\{x_{g},y_{g}\}$ of König type without squares with two $e_{i}$'s.

Finally, in section 6 we characterize the unicyclic graph with each one of the following properties: vertex decomposable, shellable, Cohen-Macaulay and well-covered.

\section{Critical, extendable and shedding vertices.}

Let $X$ be a subset of $V(G)$, the \textit{subgraph induced} by $X$ in $G$, denoted by $G[X]$ is the graph  with vertex set $X$ and whose edge set is $\{\{x,y\}\in E(G)\ |\ x,y\in X\}$. Furthermore, $G\setminus X$ denote the induced subgraph $G[V(G)\setminus X]$. Now, if $v\in V(G)$ then the set of \textit{neighbours} of $v$ is denoted by $N_{G}(v)$ and its closed neighbourhood is $N_{G}[v]=N_{G}(v) \cup \{v\}$. The \textit{degree} of $v$ in $G$ is deg$_{G}(v)=\vert N_{G}(v)\vert$. 

The contraction of the some vertex $x$ in $G$ denoted by $G/x$, is equivalent to the induced subgraph $G\setminus N_{G}[x]$. If $G$ is unmixed, shellable, sequentially Cohen-Macaulay or vertex decomposable, then $G\setminus N_{G}[x]$ is unmixed, shellable, sequentially Cohen-Macaulay or vertex decomposable, respectively. (See \cite{Biermann}, \cite{VT})

\begin{Definition}  $G$ is vertex decomposable if $G$ is a  totally disconnected graph or there is a vertex $v$ such that 

\begin{description}
 \item{\rm (a)} $G\setminus v$ and $G\setminus N_{G}[v]$ are both vertex decomposable, and
 \item{\rm (b)} each stable set in $G\setminus N_{G}[v]$ is not a maximal stable set in $G\setminus v$.
\end{description}

\end{Definition}

A \textit{shedding vertex} of $G$ is any vertex which satisfies the condition (b). Equivalently, a shedding vertex is a vertex that satisfies: For every stable set $S$ contained in $G\setminus N_{G}[v]$, there is some $x\in N_{G}(v)$ such that $S\cup \{x\}$ is stable.

\begin{Lemma} \label{4} If $x$ is a vertex of $G$, then $x$ is a shedding vertex if and only if $|N_{G}(x)\setminus N_{G}(S)|\geq 1$ for every stable set $S$ of $G\setminus N_{G}[x]$.
\end{Lemma}

\proof $\Rightarrow$) We take a stable set $S$ of $G\setminus N_{G}[x]$. Since $x$ is a shedding vertex then there is a vertex $z\in N_{G}(x)$ such that $S\cup \{z\}$ is stable set of $G\setminus x$. Thus, $z\notin N_{G}[S]$. Therefore, $|N_{G}(x)\setminus N_{G}(S)|\geq 1$. 

$\Leftarrow$) We take a stable set $S$ of $G\setminus N_{G}[x]$. Thus, there exists a vertex $z\in N_{G}(x)\setminus N_{G}(S)$. Since $z\in N_{G}(x)$, we have that $z\notin S$. Furthermore, $z\notin N_{G}(S)$ then, $S\cup \{z\}$ is a stable set of $G\setminus x$. Consequently, $S$ is not a maximal stable set of $G\setminus x$. Therefore, $x$ is a shedding vertex. \QED

\begin{Definition} Let $S$ be a stable set of $G$. If $x$ is of degree zero in $G\setminus N_{G}[S]$, then $x$ is called \textit{isolated vertex} in $G\setminus N_{G}[S]$, or we say that $S$ \textit{isolates} to $x$. \end{Definition}

By Lemma \ref{4} we have that $x$ is not a shedding vertex if and only if there exists a stable set $S$ of $G\setminus N_{G}[x]$ such that $N_{G}(x)\subseteq N_{G}(S)$, i.e. $x$ is an isolated vertex in $G\setminus N_{G}[S]$.

\begin{Corollary} Let $S$ be a stable set of $G$. If $S$ isolates to $x$ in $G$ then $x$ is not a shedding vertex in $G\setminus N_{G}[y]$ for all $y\in S$.
\end{Corollary}

\proof Since $S$ isolates to $x$ then, $S$ is a stable set in $G$ and deg$_{G\setminus N_{G}[S]}(x)=0$. Thus, $N_{G}(x)\subseteq N_{G}[S]\setminus S$Hence, if $y\in S$ and $G'=G\setminus N_{G}[y]$ then $x\in V(G')$ and $S\cap N_{G}[x]=\emptyset$. Consequently, $S'=S\setminus y$ is a stable set in $G'\setminus N_{G'}[x]$. Now, if $a\in N_{G'}(x)$ then there exists $s\in S'$ such that $\{a,s\}\in E(G')$. This implies $\vert N_{G'}(x)\setminus N_{G'}(S') \vert =0$. Therefore, by Lemma \ref{4}, $x$ is not a shedding vertex in $G'$. \QED

\begin{Theorem} \label{ext} If $x$ is a shedding vertex of $G$ then one of the following condition holds: 
\begin{description}
 \item{\rm (a)} There is $y\in N_{G}(x)$ such that $N_{G}[y]\subseteq N_{G}[x]$.
 \item{\rm (b)} $x$ is in a 5-cycle.
\end{description}
\end{Theorem}

\proof We take $N_{G}(x)=\{y_{1},y_{2},...,y_{k}\}$. If there is a vertex $y_{i}$ such that $N_{G}[y_{i}]\subseteq N_{G}[x]$ then $G$ satisfies (a). On the other hand, there are $\{z_{1},...,z_{k}\}\subseteq V(G)\setminus N_{G}[x]$ such that $\{y_{i},z_{i}\}\in E(G)$ for $i\in\{1,...,k\}$. We can suppose that $L=\{z_{1},z_{2},...,z_{q}\}=\{z_{1},...,z_{k}\}$ with $z_{i}\neq z_{j}$ for $1\leq i < j \leq q$. By Lemma \ref{4} if $L$ is a stable set of $G$ then $\vert N_{G}(x)\setminus N_{G}(L)\vert \geq 1$. But $N_{G}(x)=\{y_{1},...,y_{k}\}\subseteq N_{G}(L)$ then $L$ is not a stable set. Hence, $k\geq 2$ and there exist $z_{i_1},z_{i_2}\in L$ such that $\{z_{i_1},z_{i_2}\}\in E(G)$. Thus, there exist $y_{j_1}$ and $y_{j_2}$ such that $y_{j_1}\neq y_{j_2}$ and $\{y_{j_1},z_{i_1}\}, \{y_{j_2},z_{i_2}\}\in E(G)$. Therefore, $(x,y_{j_1},z_{i_1},z_{i_2},y_{j_2})$ is a 5-cycle of $G$.  \QED

\begin{Definition}  A vertex $v$ is called simplicial if the induced subgraph $G[N_{G}(v)]$ is a complete graph or clique. Equivalently, a simplicial vertex is a vertex that appears in exactly one clique.  \end{Definition}

\begin{Remark}\label{25} If $v,w \in V(G)$ such that $N_{G}[v]\subseteq N_{G}[w]$ then $w$ is a shedding vertex of G {\rm(see Lemma 6 in \cite{RW})}. In particular, if $v$ is a simplicial vertex then any $w\in N_{G}(v)$ is a shedding vertex {\rm(see Corollary 7 in \cite{RW})}. \end{Remark}

\begin{Corollary} Let $G$ be graph without 4-cycles. If $x$ is a shedding vertex $x$ of $G$ then $x$ is in a 5-cycle or there exists a simplicial vertex $z$ such that $\{x,z\}\in E(G)$ with $|N_{G}[z]|\leq 3$. 
\end{Corollary}

\proof By Theorem \ref{ext}, if $x$ is not in a 5-cycle then there is $z\in N_{G}(x)$ such that $N_{G}[z]\subseteq N_{G}[x]$. If deg$_{G}(z)=1$ then $z$ is a simplicial vertex. If deg$_{G}(z)=2$ then $N_{G}(z)=\{x,w\}$. Consequently, $(z,x,w)$ is a 3-cycle and $z$ is a simplicial vertex. Now, if deg$_{G}(z)\geq 3$ then there are $w_{1},w_{2}\in N_{G}(z)\setminus x$. Thus, $w_{1},w_{2} \in N_{G}(x)$ since $N_{G}[z]\subseteq N_{G}[x]$. Hence, $(w_{1},z,w_{2},x)$ is a 4-cycle of $G$. This is a contradiction. Therefore, $\vert N_{G}[z]\vert \leq 3$ and $z$ is a simplicial vertex. \QED

\begin{Remark} \label{11} If $G$ is a 5-cycle with $V(G)=\{x_{1},x_{2},x_{3},x_{4},x_{5}\}$ then $x_{i}$ is a shedding vertex. We can assume that $i=1$, then $\{x_{3}\}$ and $\{x_{4}\}$ are the stable sets in $G\setminus N_{G}[x_{1}]$. Furthermore, $\{x_{3},x_{5}\}$ and $\{x_{2},x_{4}\}$ are stable sets in $G\setminus x_{1}$. Hence, each stable set of $G\setminus N_{G}[x_1]$ is not a maximal stable set in $G\setminus x_1$. Therefore, $x_{1}$ is a shedding vertex.  \end{Remark}

\begin{Definition} A vertex $v$ of $G$ is critical if $\tau (G\setminus v)< \tau (G)$. Furthermore, $G$ is called a $vertex$ $critical$ graph if each vertex of $G$ is critical. 
\end{Definition}

\begin{Remark} \label{10} If $\tau(G\setminus v)< \tau(G)$ then $\tau(G)=\tau (G\setminus v)+1$. Moreover, $v$ is a critical vertex if and only if $\beta(G)=\beta(G\setminus v)$. \end{Remark}

\proof If $t$ is a minimal cover vertex such that $\vert t\vert=\tau(G\setminus v)$ then $t\cup v$ is a vertex cover of $G$. Thus, $\tau(G)\leq \vert t\cup v\vert =\tau(G\setminus v)+1$. But, if $\tau(G)\geq \tau(G\setminus v)+1$ then $\tau(G)=\tau(G\setminus v)+1$. 

Now, we have that $\tau (G)+\beta(G)=\vert V(G)\vert =\vert V(G\setminus v)\vert +1=\tau (G\setminus v)+\beta(G\setminus v)+1$. Hence, $\beta(G)=\beta(G\setminus v)$ if and only if $\tau (G) = \tau (G\setminus v)+1$. Therefore, $v$ is a critical vertex if and only if $\beta(G)=\beta(G\setminus v)$. \QED

\begin{Definition} A vertex $v$ of $G$ is called an $extendable$ vertex if $G$ and $G\setminus v$ are well-covered graphs such that $\beta (G)=\beta (G\setminus v)$. 
\end{Definition}

Note that if $v$ is an extendable vertex then every maximal stable set $S$ of $G\setminus v$ contains a vertex of $N_{G}(v)$.

\begin{Corollary} \label{18} Let $G$ be an unmixed graph and $x\in V(G)$. The following condition are equivalents: 

\begin{description}
 \item{\rm (a)} $x$ is an extendable vertex.
 \item{\rm (b)} $|N_{G}(x)\setminus N_{G}(S)|\geq 1$ for every stable set $S$ of $G\setminus N_{G}[x]$.
 \item{\rm (c)} $x$ is a shedding vertex.
 \item{\rm (d)} $x$ is a critical vertex and $G\setminus x$ is unmixed.
\end{description}
\end{Corollary}

\proof $(a)\Leftrightarrow (b)$ (\rm \cite{Finbow}, Lemma 2).

$(b)\Leftrightarrow (c)$ By Lemma \ref{4}.

$(a)\Leftrightarrow (d)$  Since $G$ is unmixed then by Remark \ref{10}, $x$ is extendable if and only if $x$ is a critical vertex and $G\setminus x$ is unmixed. \QED

\section{Vertex decomposable and shellable properties in graphs without 3-cycles and 5-cycles. }

\begin{Definition} A 5-cycle $C$ of $G$ is called \textit{basic} if $C$ does not contain two adjacent vertices of degree three or more in $G$. 
\end{Definition}

\begin{Lemma} \label{17} If $G$ is a graph then any vertex of degree at least 3 in a basic 5-cycle is a shedding vertex.
\end{Lemma}

\proof Let $C=(x_{1},x_{2},x_{3},x_{4},x_{5})$ be a basic 5-cycle. We suppose that deg$_{G}(x_{1})$ $\geq 3$, since $C$ is a basic 5-cycle then $\deg_{G}(x_2)=\deg_{G}(x_5)=2$. We take a stable set $S$ of $G\setminus N_{G}[x_{1}]$. Since $\{x_{3},x_{4}\}\in E(G)$ then $\vert S\cap \{x_{3},x_{4}\}\vert \leq 1$. If $x_{3}\notin S$ or $x_{4}\notin S$ then $S\cup \{x_{2}\}$ or $S\cup \{x_{5}\}$, respectively, is a stable set of $G\setminus x_{1}$. Therefore, $x_{1}$ is a shedding vertex. \QED

\begin{Definition} A simplicial complex $\Delta$ is shellable if the facets (maximal faces) of $\Delta$ can be ordered $F_{1},...,F_{s}$ such that for all $1\leq i<j\leq s$, there exists some $v\in F_{j}\setminus F_{i}$ and some $l\in \{1,...,j-1\}$ with $F_{j}\setminus F_{l}=\{v\}$. Furthermore, $F_{1},...,F_{s}$ is called  a shelling of $\Delta$. The facet set of $\Delta$ is denoted by $\mathcal{F}(\Delta)$. A graph $G$ is called shellable if $\Delta_{G}$ is shellable. \end{Definition}

\begin{Proposition} \label{26} If $G$ is shellable and $x$ is a vertex of $G$ then $G\setminus N_{G}[x]$ is shellable.
\end{Proposition}

\proof (\cite{VT}, Lemma 2.5). \QED

\begin{Lemma} \label{vec} Let $G$ be a graph such that $\{z_{1},...,z_{r}\}$ is a stable set. If $N=\bigcup_{i=1}^{r}N_{G}[z_i]$ and $G^{i+1}=G^{i}\setminus N_{G^i}[z_i]$, with $G^{1}=G$ then:
\begin{description}
 \item{\rm (a)} $\{z_{i+1},...,z_{r}\}\subseteq V(G^{i+1})$ and
 \item{\rm (b)} $G^{r+1}=G\setminus N$.
\end{description}
\end{Lemma}

\proof (a) By induction on $i$. If $i=1$, since $\{z_{1},...,z_r\}$ is a stable set then $\{z_{2},...,z_r\}\cap N_{G}[z_1]=\emptyset$ and $\{z_{2},...,z_{r}\}\subseteq V(G\setminus N_{G}[z_1])=V(G^2)$. Now, by induction hypothesis we have that $\{z_{i},z_{i+1},...,z_r\}\subseteq V(G^i)$. Since $\{z_{i},...,z_{r}\}$ is a stable set then $\{z_{i+1},...,z_{r}\}\cap N_{G^{i}}[z_{i}]=\emptyset $. Hence, $\{z_{i+1},...,z_{r}\}\subseteq V(G^{i}\setminus N_{G^{i}}[z_{i}])=V(G^{i+1})$.

(b) By induction on $r$. By $r=1$, $N=N_{G}[z_1]$ then $G^2=G\setminus N_{G}[z_1]=G\setminus N$. Now, if $r\geq 2$, we take $N^1=\bigcup_{i=1}^{r-1}N_{G}[z_i]$ and by induction hypothesis we have that $G^r=G\setminus N^1$. By incise (a), we have that $\{z_r\}\subseteq V(G^r)$. We will prove that $N^{1}\cup N_{G^{r}}[z_{r}]=N$. We have that $N^{1}\subseteq N$. Furthermore, if $y\in N_{G^{r}}[z_{r}]$ then $\{y,z_{r}\}\in E(G^{r})$. Thus, $\{y,z_{r}\}\in E(G)$ and $y\in N_{G}[z_{r}]$. Consequently, $y\in N$. Now, if $y\in N\setminus N^{1}$ then $y\in N_{G}[z_{r}]$ and $y\notin N_{G}[z_{i}]$ for $i\in \{1,...,r-1\}$. This implies that $y,z_{r}\in V(G^{r})$ and $\{y,z_{r}\}\in E(G^{r})$ implying $y\in N_{G^{r}}[z_{r}]$. Hence, $N=N^{1}\cup N_{G^{r}}[z_{r}]$. Therefore, $G^{r+1}=G^{r}\setminus N_{G^r}[z_r]=(G\setminus N^1)\setminus N_{G^r}[z_r]=G\setminus (N^{1}\cup N_{G^r}[z_r])=G\setminus N$. \QED

\begin{Definition} A subgraph $H$ of $G$ is called a c-minor (of $G$) if there exists a stable set $S$ of $G$, such that $H=G\setminus N_{G}[S]$.
\end{Definition}

\begin{Proposition} \label{30} Let $H$ be a c-minor of $G$. If $G$ is shellable then $H$ is shellable. 
\end{Proposition}

\begin{Remark} \label{35} In {\rm(\cite{wachs}, Lemma 6)} is proved that: if $G$ has a shedding vertex $v$ and if both $G \setminus v$ and $G\setminus N_{G}[v]$ are shellables with shelling $F_{1},...,F_{k}$ and $G_{1},...,G_{q}$, respectively, then $G$ is shellable. The order of the shelling is $F_{1},...,F_{k},G_{1}\cup\{v\},...,G_{q}\cup\{v\}$. \end{Remark}

\begin{Theorem} \label{8} Let $G$ be a connected graph with a basic 5-cycle $C$. $G$ is a shellable graph if and only if there is a shedding vertex $x\in V(C)$ such that $G\setminus x$ and $G\setminus N_{G}[x]$ are shellable graphs.
\end{Theorem}
 
\proof $\Rightarrow$) We can suppose that $C=(x_{1},x_{2},x_{3},x_{4},x_{5})$. If $G=C$ then $G$ is shellable. By Remark \ref{11} each vertex is a shedding vertex. Furthermore, $G_{1}=G\setminus x_{1}$ is a path with shelling $\{x_{2},x_{4}\}, \{x_{2},x_{5}\}, \{x_{3},x_{5}\}$ and $G\setminus N_{G}[x_{1}]$ is an edge. Therefore, $G_{1}$ and $G\setminus N_{G}[x_{1}]$ are shellable graphs. Now, we suppose $G\neq C$. We can assume that deg$_{G}(x_{1})\geq 3$ then, since $C$ is a basic 5-cycle deg$_{G}(x_{2})=$ deg$_{G}(x_{5})=2$. Also, we can suppose deg$_{G}(x_{3})=2$ and deg$_{G}(x_{4})\geq 2$. By Lemma \ref{17}, $x_{1}$ is a shedding vertex. $G\setminus N_{G}[x_{1}]$ is a shellable graph by Proposition \ref{26}. We will prove that $G_{1}$ is shellable. Since $G$ is shellable then $G_{2}=G\setminus N_{G}[x_{2}]$ is shellable. We assume that $F_{1},...,F_{r}$ is a shelling of $\Delta_{G_{2}}$. Now, we take $G_{3}=G\setminus N_{G}[x_{3},x_{5}]$ then $G_{3}$ is shellable and whose shelling is $H_{1}, H_{2},...,H_{k}$. We take $F\in \mathcal{F}(\Delta_{G_1})$. If $x_{2}\in F$ then $F\setminus x_{2}\in \mathcal{F}(\Delta_{G_{2}})$ and there exists $F_{i}$ such that $F=F_{i}\cup x_{2}$. If $x_{2}\notin F$ then $x_{3}\in F$ and $x_{4}\notin F$. Thus, $x_{5}\in F$. Hence, $F\setminus \{x_{3},x_{5}\}\in \mathcal{F}(\Delta_{G_{3}})$ then there exist $H_{j}$ such that $F=H_{j}\cup \{x_{3},x_{5}\}$. This implies, $\mathcal{F}(\Delta_{G_{1}})=\{F_{1}\cup \{x_{2}\},...,F_{r}\cup \{x_{2}\}, H_{1}\cup \{x_{3},x_{5}\},...,H_{k}\cup \{x_{3},x_{5}\}\}$. Furthermore, $F_{1}\cup \{x_{2}\},...,F_{r}\cup \{x_{2}\}$ and $H_{1}\cup \{x_{3},x_{5}\},...,H_{k}\cup \{x_{3},x_{5}\}$ are shellings. Now, $x_{3}\in (H_{j}\cup \{x_{3},x_{5}\}) \setminus (F_{i}\cup \{x_{2}\})$ and $H_{j}$ is a stable set of $G$ without vertices of $C$ then $H_{j}\cup \{x_{2},x_{5}\}$ is a maximal stable set of $G$ since  $N_{G}(x_{2},x_{5})=V(C)$ and $\{x_{2},x_{5}\}\notin E(G)$. Consequently, $H_{j}\cup \{x_{2},x_{5}\}=F_{l}\cup \{x_{2}\}$ for some $l\in \{1,...,r\}$ then $(H_{j}\cup \{x_{3},x_{5}\})\setminus (F_{l}\cup \{x_{2}\})=\{x_{3}\}$. Therefore, $G_{1}$ is a shellable graph. 

$\Leftarrow$) By Remark \ref{35}. \QED

\begin{Definition} If $v,w\in V(G)$ then the distance d$(u,v)$ between $u$ and $v$ in $G$ is the length of the shortest path joining them, otherwise d$(u,v)=\infty$. Now, if $H\subseteq G$ then the distance from a vertex $v$ to $H$ is d$(v,H)=min\{d(v,u)\ \vert \ u\in V(H)\}$. Furthermore, we define $D_{i}(H)=\{v\in V(G)\ \vert \ d(v,H)=i\}$. \end{Definition}

\begin{Lemma} Let $G$ be a graph. If  $W \subseteq V(G)$ then $D_0(W)=W$ and $D_i(W)=\{x \in V(G)\ |\  x \notin D{_j}(W) \mbox{ for } 0 \leq j < i \mbox{ and } x\in N_G(D_{i-1}(W)) \}$ for $i>0$. \end{Lemma}

\proof We take $D_{i}=D_{i}(W)$. By definition $D_{0}=\{v\in V(G)\ \vert \ d(v,W)=0\}=W$. Now, we take $A=\{x \in V(G)\ |\  x \notin D_j \mbox{ for } 0 \leq j < i \mbox{ and } x\in N_G(D_{i-1}) \}$. If $v\in D_{i}$ then d$(v,W)=i$. Consequently, there is a vertex $z$ such that $\{v,z\}\in E(G)$ and d$(z,W)=i-1$. Hence, $v\in N_{G}(D_{i-1})\setminus N_{G}(D_{j})$ for $0\leq j<i$. Therefore, $v\in A$. On the other hand, $v\in A$ then d$(v,W)\geq i$ since $v\notin D_{j}$ for $0\leq j<i$ and d$(v,W)\leq i$ since $v\in N_{G}(D_{i-1})$. Therefore, d$(v,W)=i$. \QED

\begin{Definition} A \textit{cut vertex} of a graph is one whose removal increases the number of connected components. A \textit{block} of a graph is a maximal subgraph without cut vertices. A connected graph without cut vertices with at least three vertices is called 2-connected graph.  \end{Definition}

\begin{Remark}\label{19} Each connected component of a graph $G$ is a c-minor of $G$. \end{Remark}

\begin{Definition} A vertex of degree one is a \textit{leaf} or \textit{free vertex}. Furthermore, an edge which is incident with a leaf is called \textit{pendant}. 
\end{Definition}

In the following result $P$ is a property closed under c-minors, i.e., if $G$ has the property $P$ then each c-minor has the property $P$.

\begin{Theorem} \label{5} Let $G$ be a graph without 3-cycles and 5-cycles with a $2$-connected block $B$. If $G$ satisfies the property $P$ and $B$ does not satisfy it then there exists $x\in D_1(B)$ such that \rm deg$_G(x)=1$ 
\end{Theorem}

\proof By contradiction, we assume that if $x\in D_{1}(B)$, then $\vert N_G(x) \vert >1$. Thus, there exist $a,b \in N_G(x)$ with $a\neq b$. We can suppose that $a\in V(B)$. If $b \in V(B)$ then $G[\{x\}\cup V(B)]$ is $2$-connected, but $B\subsetneq G[\{x\}\cup V(B)]$. This is a contradiction since $B$ is a block. Consequently, $V(B)\cap N_{G}(x)=\{a\}$. Now, we suppose that $b \in D_1(B)$. Since there is not a 3-cycle in $G$, then $a\notin N_{G}(b)$. Hence, there exists $c \in N_G(b) \cap V(B)$ such that $c \neq a$. This implies $G[\{x,b\}\cup V(B)]$ is $2$-connected. But $B\subsetneq G[\{x,b\}\cup V(B)]$, a contradiction then, $D_1(B)\cap N_G(x)=\emptyset$. Thus, $N_G(x)\cap (V(B) \cup D_1(B))=\{a\}$ and $b\in D_{2}(B)$. Now, if $D_1(B)=\{x_1,\ldots,x_r\}$ then, there exist $a_i$ such that $V(B)\cap N_G(x_i)=\{a_i\}$. Furthermore, there exists $b_i$ such that  $b_i \in N_G(x_i)\cap D_2(B)$. We can suppose that $L=\{b_1,\ldots,b_r\}=\{b_1,\ldots,b_s\}$ with $b_i\neq b_j$ for $1\leq i<j\leq s$. We will prove that $L$ is a stable set. Suppose that $\{b_i, b_j\}\in E(G)$, if $a_i=a_j$ then $(a_i,x_i,b_i,b_j,x_j,a_i)$ is a 5-cycle in $G$.
This is a contradiction, consequently $a_i \neq a_j$ and the induced subgraph $G[\{x_i,b_i,b_j,x_j\}\cup V(B)]$ is $2$-connected. But $B$ is a block then,  $\{b_i, b_j\}\notin E(G)$. Therefore, $L$ is a stable set and $G'=G\setminus(N_G[L])$ satisfies the property $P$. Furthermore $D_1(B)\subset N_G(L)$ then, $B$ is a connected component of $G'$. But, $B$ does not satisfy $P$, this is a contradiction. Therefore, there exists a free vertex in $D_1 (B)$. \QED

\begin{Remark} \label{27} The following properties: unmixed, shellable, Cohen-Macaulay, sequentially Cohen-Macaulay and vertex decomposable are closed under c-minors. \rm (See \cite{Biermann}, \cite{Villa}).
\end{Remark}

\begin{Corollary} \label{3} Let $G$ be a graph without 3-cycles and 5-cycles and $B$ a $2$-connected block. If $G$ is shellable (unmixed, Cohen-Macaulay, sequentially Cohen-Macaulay or vertex decomposable) and  $B$ is not shellable (unmixed, Cohen - Macaulay, sequentially Cohen-Macaulay or vertex decomposable) then there exists $x\in D_1(B)$ such that \rm deg$_G(x)=1$ 
\end{Corollary}

\proof From Remark \ref{27} and Proposition \ref{5} there exists $x\in D_{1}(B)$ such that deg$_{G}(x)=1$.   \QED

\begin{Corollary} \label{20} Let $G$ be a bipartite graph and $B$ a 2-connected block. If $G$ is shellable then there exists $x\in D_{1}(B)$ such that $\deg_{G}(x)=1$.  \end{Corollary}

\proof Since $G$ is bipartite then $B$ is bipartite. If $H$ is a shellable bipartite graph then $H$ has a free vertex. (see \cite{VT}, Lemma 2.8). In particular, $H$ is not 2-connected. Hence, $B$ is not shellable. Therefore, by Corollary \ref{3}, there exists $x\in D_{1}(B)$ such that deg$_{G}(x)=1$. \QED

\begin{Lemma} \label{15} Let $G$  be a graph without 3-cycles and 5-cycles. If $G$ is vertex decomposable then $G$ has a free vertex.
\end{Lemma}

\proof Since $G$ is vertex decomposable then, there is a shedding vertex $x$. Since there are not 5-cycles in $G$ then by Theorem \ref{ext}, there exists $y\in N_{G}[x]$ such that $N_{G}[y]\subseteq N_{G}[x]$. If $z\in N_{G}(y)\setminus x$ then $(x,y,z)$ is a 3-cycle. This is a contradiction. Therefore, $N_{G}(y)=\{x\}$ and $y$ is a free vertex. \QED

\begin{Remark} \label{28} If $G$ is a vertex decomposable graph, then $G\setminus N_{G}[S]$ is a vertex decomposable graph for any $S$ stable set of $G$
\end{Remark} 

\proof See Theorem 2.5 in \cite{Biermann}. \QED

\begin{Theorem} Let $G$  be a graph without 3-cycles and 5-cycles and $\{x,y\}\in E(G)$ such that $x$ is a free vertex. If $G_{1}=G\setminus N_{G}[x]$ and $G_{2}=G\setminus N_{G}[y]$, then $G$ is vertex decomposable if and only if $G_{1}$ and $G_{2}$ are vertex decomposable.  
\end{Theorem} 

\proof $\Rightarrow$) By Remark \ref{28} $G_{1}$ and $G_{2}$ are vertex decomposable.  

$\Leftarrow$) Since $G_{1}$ is vertex decomposable then $G_{1}\cup \{x\}=G\setminus y$ is vertex decomposable. Furthermore, $y$ is a shedding vertex by -Remark \ref{25} and $G_{2}$ is vertex decomposable. Therefore, $G$ is vertex decomposable. \QED

\begin{Corollary} \label{14} If $G$ is a 2-connected graph without 3-cycles and 5-cycles then, $G$ is not a vertex decomposable.
\end{Corollary}

\proof If $G$ is a vertex decomposable graph then, by Lemma \ref{15}, $G$ has a free vertex. This is a contradiction since $G$ is 2-connected. Therefore, $G$ is not vertex decomposable. \QED

\begin{Definition} Let $G_{1},G_{2}$ be graphs. If $K=G_{1}\cap G_{2}$ is a complete graph with $|V(K)|=k$ then $G=G_{1}\cup G_{2}$ is called the \textit{k-clique-sum} (or clique-sum) of $G_{1}$ and $G_{2}$ in $K$.
\end{Definition}

\begin{Remark} \label{29} In {\rm(\cite{Chris}, Proposition 4.1)} is proved that if $C_{n}$ is a n-cycle then $C_{n}$ is vertex decomposable, shellable or sequentially Cohen-Macaulay if and only if $n=3$ or $5$. Furthermore, a chordal graph is vertex decomposable. {\rm(see \cite{RW}, Corollary 7 (2))} \end{Remark}

\begin{Corollary} If $G$ is the 2-clique-sum of the cycles $C_{1}$ and $C_{2}$ with $\vert V(C_{1}) \vert = r_{1}\leq r_{2}=\vert V(C_{2})\vert$ then, $G$ is vertex decomposable if and only if $r_{1}=3$ or $r_{1}=r_{2}=5$.
\end{Corollary}

\proof $\Leftarrow$) First, we suppose that $r_{1}=3$. Consequently, we can assume $C_{1}=(x_{1},x_{2},x_{3})$ and $x_{2},x_{3}\in V(C_{1})\cap V(C_{2})$. Thus, $x_{1}$ is a simplicial vertex. Hence, by Remark \ref{25} $x_{2}$ is a shedding vertex. Furthermore $G\setminus x_{2}$ and $G\setminus N_{G}[x_{2}]$ are trees then, $G\setminus x_{2}$ and $G\setminus N_{G}[x_{2}]$ are vertex decomposable graphs by Remark \ref{29}. Hence, $G$ is vertex decomposable. 

Now, we assume that $r_{1}=r_{2}=5$ with $C_{1}=(x_{1},x_{2},x_{3},x_{4},x_{5})$ and $C_{2}=(y_{1},x_{2},x_{3},y_{4},y_{5})$. We take a stable set $S$ in $G\setminus N_{G}[x_{5}]$. If $x_{2}\in S$ then, $S\cup \{x_{4}\}$ is a stable set in $G\setminus x_{5}$. If  $x_{2}\notin S$ then,  $S\cup \{x_{1}\}$ is a stable set in $G\setminus x_{5}$. Consequently, $x_{5}$ is a shedding vertex.  Since $x_{2}$ is a neighbour of a free vertex in $G_{1}=G\setminus x_{5}$ then, $x_{2}$ is a shedding vertex in $G_{1}$. Furthermore, since $G_{1}\setminus x_{2}$ and $G_{1}\setminus N_{G_{1}}[x_{2}]$ are forest then, they are vertex decomposable graphs by Remark \ref{29}. Thus, $G_{1}$ is vertex decomposable. Now, $G\setminus N_{G}[x_{5}]=C_{2}$ then it is vertex decomposable by Remark \ref{29}. Therefore, $G$ is vertex decomposable.

$\Rightarrow$) By Corollary \ref{14} we have that $\{r_{1},r_{5}\}\cap \{3,5\}\neq\emptyset$. If $r_{1}\neq 3$ then $r_{1}=5$ or $r_{2}=5$. Consequently, we can assume that  $\{C_{1},C_{2}\}=\{C,C'\}$ where $C=(x_{1},x_{2},x_{3},x_{4},x_{5})$ and $x_{2},x_{3}\in V(C)\cap V(C')$. Thus, $G\setminus N_{G}[x_{5}]=C'$ is vertex decomposable by Remark \ref{28} then $\vert V(C')\vert\in \{3,5\} $. But, $r_{1}\neq 3$  then $\vert V(C')\vert=5$ and $r_{1}=r_{5}=5$. Therefore, $r_{1}=3$ or $r_{1}=r_{2}=5$. \QED

\begin{Definition} The $girth$ of $G$ is the length of the smallest cycle or if $G$ is a forest we consider the girth of $G$ infinite. \end{Definition}

\begin{Lemma} \label{16} Let $G$ be a 2-connected graph with girth at least 11. Then $G$ is not shellable. 
\end{Lemma}

\proof Let $r$ be the girth of $G$, then there exists a cycle $C=(x_{1},x_{2},...,x_{r})$. If $G=C$ then, $G$ is not shellable. Hence, $G\neq C$ then, $D_{1}(C)\neq \emptyset$. We take $y\in D_{1}(C)$, without loss of generality we can assume that $\{x_{1},y\}\in E(G)$. If $\{x_{i},y\}\in E(G)$ for some $i\in \{2,...,r\}$ then we take the cycles $C_{1}=(y,x_{1},x_{2},...,x_{i})$ and $C_{2}=(y,x_{1},x_{r},x_{r-1},...,x_{i})$. Thus, $\vert V(C_{1})\vert =i+1$ and $\vert V(C_{2})\vert = r-i+3$ Since $r$ is the girth of $G$ then $i+1\geq r$ and $r-i+3\geq r$. Consequently, $3\geq i$ implies $4\geq r$. But $r\geq 11$, this is a contradiction. This implies that $\vert N_{G}(y)\cap V(C)\vert = 1$. Now, we suppose that there exist $y_{1},y_{2}\in D_{1}(C)$ such that $\{y_{1},y_{2}\}\in E(G)$. We can assume that $\{x_{1},y_{1}\}, \{x_{i},y_{2}\}\in E(G)$. Since $r\geq 11$ then, there are not 3-cycles in $G$. In particular, $x_{1}\neq x_{i}$. Now, we take the cycle $C'=(y_{1},x_{1},...,x_{i},y_{2})$ and $C''=(y_{1},x_{1},x_{r},x_{r-1},...,x_{i},y_{2})$. But, $\vert V(C')\vert =i+2$ and $\vert V(C'')\vert = r-i+4$. Since $r$ is the girth then $i+2\geq r$ and $r-i+4\geq r$. Hence, $4\geq i$ and $6\geq r$. This is a contradiction then $D_{1}(C)$ is a stable set. Now, since $G$ is 2-connected then for each $y\in D_{1}(C)$ there exists $z\in N_{G}(y)\cap D_{2}(C)$. Now, if there exist $z_{1},z_{i}\in D_{2}(C)$ such that $\{z_{1},z_{i}\}\in E(G)$ then, there exist $y_{1},y_{i}\in D_{1}(C)$ such that $\{z_{1},y_{1}\}, \{z_{i},y_{i}\}\in E(G)$. Since there are not 3-cycles then $y_{1}\neq y_{i}$. We can assume that $\{x_{1},y_{1}\}, \{x_{i},y_{i}\}\in E(G)$. Since there are not 5-cycles then $i\neq 1$. Consequently, there exist cycles $C_{1}'=(x_{1},...,x_{i},y_{i},z_{i},z_{1},y_{1})$ and $C_{2}'=(x_{i},...,x_{r},x_{1},y_{1},z_{1},z_{i},y_{i})$. This implies $r\leq \vert V(C_{1}')\vert =i+4$ and $r\leq \vert V(C_{2}')\vert =r-i+6$. Hence, $i\leq 6$ implies $r\leq 10$. This is a contradiction then $D_{2}(C)$ is a stable set. Thus, $C$ is a connected component of $G\setminus N_{G}[D_{2}(C)]$. But $C$ is not shellable, therefore $G$ is not shellable. \QED

\begin{Theorem} If $G$ has girth at least 11 then $G$ is shellable if and only if there exists $x\in V(G)$ with $N_{G}(x)=\{y\}$ such that $G\setminus N_{G}[x]$ and $G\setminus N_{G}[y]$ are shellables.    
\end{Theorem}

\proof $\Leftarrow$) By (\cite{VT}, Theorem 2.9).

$\Rightarrow$) Since shellability is closed under c-minors, it is only necessary to prove that there exist $x\in V(G)$, such that deg$_{G}(x)=1$. If every block of $G$ is an edge then $G$ is a forest and there exist $x\in V(G)$ with deg$_{G}(x)=1$. Hence, we can assume that there exist $B$ a 2-connected block of $G$. Since $B$ is an induced subgraph of $G$, then its girth is at least 11. Thus, by Lemma \ref{16}, $B$  is not shellable. Consequently, by Theorem \ref{5}, there exist $x\in D_{1}(B)$ such that $\deg_{G}(x)=1$. \QED

\section{König graph and well-covered graphs without 3-cycles and 5-cycles.}

In \cite{half} and \cite{MM} proved that Cohen-Macaulay, pure shellable, unmixed vertex decomposable are equivalent properties if $G$ is a very well-covered graph. Furthermore, they gave a theorem similar to theorem given in \cite{HH} by Herzog and Hibi for bipartite graphs. 

In this paper we denoted by $Z_{G}$ the set of the isolated vertices, that is, $$Z_{G}=\{x\in V(G)\ \vert \ \deg_{G}(x)=0\}.$$

\begin{Definition} A graph $G$ is called \textit{very well-covered} if it is well-covered without isolated vertices and with 2ht$(I(G))=|V(G)|$. \end{Definition}

\begin{Definition} $G$ is a \textit{König} graph if $\tau(G)=\nu(G)$ where $\nu(G)$ is the maximum number of pairwise disjoint edges. A \textit{perfect matching of König type} of $G$ is a collection $e_{1},...,e_{g}$ of pairwise disjoint edges whose union is $V(G)$ and $\tau(G)=g$. \end{Definition}

\begin{Proposition} \label{22} Let $G$ be a König graph and $G'=G\setminus Z_{G}$. Then the following are equivalent:
\begin{description}
 \item{\rm (a)} $G$ is unmixed
 \item{\rm (b)} If $G'\neq \emptyset$, then $G'$ has a perfect matching $e_{1},...,e_{g}$ of König type such that for any two edges $f_{1}\neq f_{2}$ and for two distinct vertices $x\in f_{1}$, $y\in f_{2}$ contained in some $e_{i}$, one has that $(f_{1}\setminus x)\cup (f_{2}\setminus y)$ is an edge.
\end{description}
\end{Proposition}

\proof By (\cite{MR}, Lemma 2.3 and Proposition 2.9). \QED

\begin{Lemma} \label{12} $G$ is an unmixed König graph if and only if $G$ is totally disconnected or $G^{'}=G\setminus Z_{G}$ is very well-covered.
\end{Lemma}

\proof $\Rightarrow$) If $G$ is not totally disconnected, then from Proposition \ref{22} $G'$ has a perfect matching $e_{1},...,e_{g}$ of König type. Hence, $|V(G')|=2\tau(G')=2ht(I(G'))$. Furthermore, $G'$ is unmixed then $G'$ is well-covered.

$\Leftarrow$) If $G$ is totally disconnected then $\nu(G)=0$ and $\tau(G)=0$. Hence, $G$ is unmixed König graph. Now, if $G$ is not totally disconnected then $G'$ is very well-covered. Consequently, by (\cite{GV}, Remark 2.2) $G'$ has a perfect matching. Thus, $\nu (G')= \vert V(G')\vert /2=ht (G')=\tau(G')$. Furthermore, $G'$ is well-covered. Hence, $G'$ is unmixed König graph. Therefore, $G$ also is unmixed König graph. \QED

\begin{Lemma} \label{P3} If $G$ is a well-covered graph without 3-cycles, 5-cycles and 7-cycles, then $G$ is a König graph.
\end{Lemma}

\proof By induction on $\vert V(G)\vert$. We take $x\in V(G)$ then, by Remark \ref{27} $G_{1}=G\setminus N_{G}[x]$ is an unmixed graph. Furthermore, $G_{1}$ does not contain 3-cycles, 5-cycles and 7-cycles then, by induction hypothesis, $G_{1}$ is König, by Proposition \ref{22}. If $G_{2}=G_{1}\setminus Z_{G}=\emptyset$ then $V(G)=N_{G}[x]\cup Z_[G_{1}]$. Hence, $x\cup Z_{G_{1}}$ and $N_{G}(x)$ are stable set and $G$ is bipartite. Consequently, $G$ is König. Therefore, we can assume that $G_{2}\neq \emptyset$ implying $G_{2}$ has a perfect matching $e_{1}=\{x_{1},y_{1}\},...,e_{g}=\{x_{g},y_{g}\}$ of König type. We can assume that $N_{G}(x)=\{z_{1},...,z_{r}\}$ and $D=\{x_{1},...,x_{g}\}$ is a minimal vertex cover of $G$. Thus, $F=\{y_{1},...,y_{g}\}$ is a maximal stable set of $G_{2}$. Now, we take the subsets: $A_{1}=N_{G}(z_{1},...,z_{r})\cap D$, $B_{1}=\{y_{j}\in F\ \vert\ x_{j}\in A_{1}  \}$, $B_{2}=N_{G}(z_{1},...,z_{r})\cap F$ and $A_{2}=\{x_{j}\in D\ \vert\ y_{j}\in B_{2}  \}$. If there exists $y_{i}\in B_{1}\cap B_{2}$ then $x_{i}\in A_{1}$ and there exist $z_{k},z_{p}\in N_{G}(x)$ such that $\{x_{i},z_{k}\}, \{y_{i},z_{p}\}\in E(G)$. If $k=p$ then $(z_{k},x_{i},y_{i})$ is a 3-cycle and if $k\neq p$ then $(x,z_{k},x_{i},y_{i},z_{p})$ is a 5-cycle. This is a contradiction, consequently $B_{1}\cap B_{2}=\emptyset$. Now, we take the subsets $B_{3}=(N_{G}(A_{2})\cap F)\setminus B_{2}$, $A_{3}=\{x_{j}\in D\ \vert\ y_{j}\in B_{3}  \}$, $B_{4}=(N_{G}(A_{1})\cap F)\setminus B_{1}$ and $A_{4}=\{x_{j}\in D\ \vert\ y_{j}\in B_{4}  \}$. If $y_{i}\in B_{1}\cap B_{3}$ then there exist $x_{j}\in A_{2}$ and $z_{k},z_{p}\in N_{G}(x)$ such that $\{x_{i},z_{k}\}, \{y_{i},x_{j}\}$ and $\{y_{j},z_{p}\}\in E(G)$. Hence, if $k=p$ then $(z_{k},x_{i},y_{i},x_{j},y_{j})$ is a 5-cycle and if $k\neq p$ then $(x,z_{k},x_{i},y_{i},x_{j},y_{j},z_{p})$ is a 7-cycle. A contradiction, implying $B_{1}\cap B_{3}=\emptyset$. Now, if $y_{i}\in B_{2}\cap B_{4}$ then there exist $x_{j}\in A_{1}$ and $z_{k},z_{p}\in N_{G}(x)$ such that $\{x_{j},z_{k}\}, \{x_{j},y_{i}\}$ and  $\{y_{i},z_{p}\}\in E(G)$. Consequently, if $k=p$ then $(z_{k},x_{j},y_{i})$ is a 3-cycle and if $k\neq p$ then $(x,z_{k},x_{j},y_{i},z_{p})$ is a 5-cycle. This is a contradiction, hence $B_{2}\cap B_{4}=\emptyset$. Now, if $y_{i}\in B_{3}\cap B_{4}$ then there exist $x_{j}\in A_{1}$, $x_{q}\in A_{2}$ and $z_{k},z_{p}\in N_{G}(x)$ such that $\{x_{j},y_{i}\}$, $\{x_{q},y_{i}\}$, $\{x_{j},z_{k}\}$ and $\{y_{q},z_{p}\}\in E(G)$. Thus, if $k=p$ then $(z_{k},x_{j},y_{i},x_{q},y_{q})$ is a 5-cycle and if $k\neq p$ then $(x,z_{k},x_{j},y_{i},x_{q},y_{q},z_{p})$ is a 7-cycle. This is a contradiction, implying $B_{3}\cap B_{4}=\emptyset$. Therefore, $B_{1},B_{2},B_{3},B_{4},B_{5}$ are pairwise disjoint sets, where $B_{5}=F\setminus (B_{1}\cup B_{2}\cup B_{3}\cup B_{4})$. 

Now, we will prove that if $A'=A_{1}\cup A_{4}\cup A_{5}$ and $B'=B_{1}\cup B_{4}\cup B_{5}$ then $N_{G}(B')\subseteq A'$. First, we observe that $N_{G}(Z_{G_{1}})\subseteq N_{G}(x)$. Since $B_{3}\cap B_{5}=\emptyset$ then $N_{G}(B_{5})\cap A_{2}=\emptyset$. Furthermore, if $x_{i}\in N_{G}(B_{5})\cap A_{3}$ then, there exist $y_{j}\in B_{5}$ and $x_{q}\in A_{2}$ such that $\{x_{i},y_{j}\}$, $\{x_{q},y_{i}\}\in E(G)$. Hence, since $G_{1}$ is unmixed and by Proposition \ref{22} $\{x_{q},y_{j}\}=(\{x_{q},y_{i}\}\setminus y_{i})\cup(\{x_{i},y_{j}\}\setminus x_{i})\in E(G)$. Consequently, $y_{j}\in B_{3}\cap B_{5}$, this is a contradiction. Therefore,  $N_{G}(B_{5})\cap A_{3}=\emptyset$ and $N_{G}(B_{5})\subseteq A'$. Also, since $B_{1}\cap B_{3}=\emptyset$ then $N_{G}(B_{1})\cap A_{2}=\emptyset$. Furthermore, if $x_{i}\in N_{G}(B_{1})\cap A_{3}$ then there exist $y_{j}\in B_{1}$, $x_{q}\in A_{2}$ such that $\{x_{i},y_{j}\}, \{x_{q},y_{i}\}\in E(G)$. Since $G_{1}$ is unmixed, by Proposition \ref{22} $\{x_{q},y_{j}\}=(\{x_{i},y_{j}\}\setminus x_{i}) \cup (\{x_{q},y_{i}\}\setminus y_{i})\in E(G)$. Then $y_{j}\in B_{1}\cap B_{3}$. This is a contradiction, thus $N_{G}(B_{1})\cap A_{3}=\emptyset$, and $N_{G}(B_{1})\subseteq A'$. Now, since $B_{3}\cap B_{4}=\emptyset$ then $N_{G}(B_{4})\cap A_{2}=\emptyset$. Furthermore, if $x_{i}\in N_{G}(B_{4})\cap A_{3}$ then there exist $y_{j}\in B_{4}$ and $x_{q}\in A_{2}$ such that $\{x_{i},y_{j}\}$ and $\{x_{q},y_{i}\}\in E(G)$. Consequently, since $G_{1}$ is unmixed $\{x_{q},y_{j}\}=(\{x_{i},y_{j}\}\setminus x_{i}) \cup (\{x_{q},y_{i}\}\setminus y_{i})\in E(G)$. This is a contradiction since $B_{3}\cap B_{4}=\emptyset$. Hence, $N_{G}(B_{4})\cap A_{3}=\emptyset$ and $N_{G}(B_{4})\subseteq A'$. Therefore, $N_{G}(B')=A'$.

This implies $G_{3}=G\setminus N_{G}[B']=G\setminus (A'\cup B')$. Furthermore, since $B'$ is a stable set then, $G_{3}$ is unmixed without 3-cycles, 5-cycles and 7-cycles. Now, if $G\neq G_{3}$ then, by induction hypothesis $G_{3}$ is König. That is, $\tau(G_{3})=\nu(G_{3})$. Consequently, if $D'$ is a minimal vertex cover of $G_{3}$ then $D'\cup A'$ is a vertex cover of $G$ since $N_{G}(B')=A'$. Furthermore, $G[A'\cup B']$ has a perfect matching with $\vert A' \vert$ elements. Thus, $\tau(G) \leq \tau(G_{3})+\vert A' \vert =\nu(G_{3})+\vert A' \vert\leq \nu(G)$. Hence, $\tau(G)=\nu(G)$, implying, $G$ is König. Therefore, we can assume $G=G_{3}$.

Now, if there exist $x_{i},x_{j}\in A_{2}$ such that $\{x_{i},x_{j}\}\in E(G)$ then there exist $z_{k},z_{p}\in N_{G}(x)$ such that $\{y_{i},z_{k}\}$ and $\{y_{j},z_{p}\}\in E(G)$. If $k=p$ then $(z_{k},y_{i},x_{i},x_{j},y_{j})$ is a 5-cycle and if $k\neq p$ then $(x,z_{k},y_{i},x_{i},x_{j},y_{j},z_{p})$ is a 7-cycle. This is a contradiction, then $A_{2}$ is a stable set. Similarity, if  $\{x_{i},x_{j}\}\in E(G)$ with $x_{i}\in A_{2}$ and $x_{j}\in A_{3}$, then there exist $x_{q}\in A_{2}$ such that $\{x_{q},y_{j}\}\in E(G)$. If $q=i$ then $(x_{i},x_{j},y_{j})$ is a 3-cycle. This is a contradiction, hence $q\neq i$. Consequently, since $G_{3}$ is unmixed and by Proposition \ref{22}, $\{x_{i},x_{q}\}=(\{x_{i},x_{j}\}\setminus x_{j})\cup(\{x_{q},y_{j}\}\setminus y_{j})\in E(G_{3})$. This is a contradiction since $A_{2}$ is a stable set of $G$. Thus, there are not edges between $A_{2}$ and $A_{3}$. Now, $\{x_{i},x_{j}\}\in E(G)$ with $x_{i},x_{j}\in A_{3}$ then there is a vertex $x_{q}\in A_{2}$ such that $\{x_{q},y_{i}\}\in E(G)$. Hence, since $G_{3}$ is unmixed then $\{x_{j},x_{q}\}=(\{x_{i},x_{j}\}\setminus x_{i}) \cup (\{x_{q},y_{i}\}\setminus y_{i})\in E(G_{3})$. But, there are not edges between $A_{2}$ and $A_{3}$. Therefore, $A_{3}$ is a stable set, implying $A_{2}\cup A_{3}$ is also a stable set in $G$. This implies, $\{x\}\cup Z_{G_{1}}\cup B_{2}\cup B_{3}$ and $N_{G}(x)\cup A_{2}\cup A_{3}$ are stable sets, since $N_{G}(Z_{G_{1}})\subseteq N_{G}(x)$, implying $G$ is bipartite. Therefore, $G$ is König. \QED

\begin{Theorem} \label{24} If $G$ is a graph without 3-cycles, 5-cycles and 7-cycles then the following conditions are equivalent:
\begin{description}
 \item{\rm (1)} $G$ is well-covered.
 \item{\rm (2)} $G$ has a perfect matching $e_{1},...,e_{g}$ of König type such that if $f_{1},f_{2}\in E(G)$ with $a\in f_{1}$, $b\in f_{2}$ and $\{a,b\}=e_{i}$ then, $(f_{1}\setminus a)\cup (f_{2}\setminus b)\in E(G)$. 
\end{description}
\end{Theorem}

\proof (1) $\Rightarrow$ (2) By Lemma \ref{P3}, $G$ is König. Furthermore, since $G$ is well-covered then by Proposition \ref{22} $G$ satisfies (2).

(2) $\Rightarrow$ (1) By Proposition \ref{22}, $G$ is well covered. \QED

\begin{Example}
The subgraph $G$ is unmixed without 3-cycles and 5-cycles. If $G$ is König then $G$ has a perfect matching by {\rm (\cite{MR}, Lemma 2.3)} and $\vert V(G)\vert $ is even. But $\vert V(G)\vert=9$, thus $G$ is not König.  

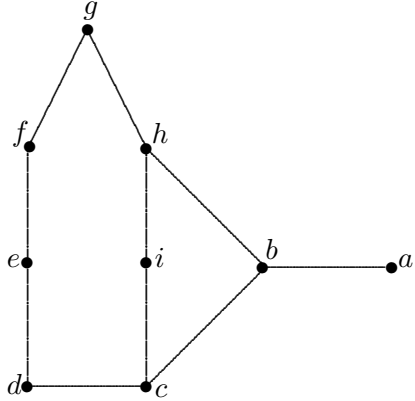
\begin{center}

\begin{picture} (300,100)
\plot 0 -45 0 45 22.5 90 45 45 45 -45 0 -45 /
\plot  45 45 90 0 45 -45  /
\plot 90 0 135 0 /

\put(-2,43){$\bullet$}   
\put(-3,-1){$\bullet$} 
\put(42,42){$\bullet$}  
\put(42,-1){$\bullet$} 
\put(42,-48){$\bullet$}  
\put(86,-3){$\bullet$} 
\put(-3,-48){$\bullet$}  
\put(135,-3){$\bullet$} 
\put(20,87){$\bullet$}

\put(-6,48){$f$}  
\put(-8,0){$e$} 
\put(47,47){$h$}  
\put(48,-49){$c$}  
\put(90,3){$b$} 
\put(-8,-48){$d$}  
\put(140,0){$a$} 
\put(48,0){$i$} 
\put(21,96){$g$}

\end{picture}
\end{center}

\end{Example} 

\vspace{18mm}

\begin{Lemma} \label{21} If $G$ is unmixed graph then for each $x\in V(G)$, $N_{G}(x)$ does not contain two free vertices.
\end{Lemma}

\proof By contradiction suppose that there exists $x\in V(G)$ such that $y_{1},...,y_{s}$ are free vertices in $N_{G}(x_{1})$ with $s\geq 2$. Hence, $G'=G\setminus N_{G}[y_{1},...,y_{s}]=G\setminus \{x,y_{1},...,y_{s}\}$ is unmixed. Now, we take $S$ a maximal stable set of $G'$. Thus, $\vert S\vert=\beta(G')$ since $G'$ is unmixed. Consequently, $S_{1}=s\cup \{y_{1},...,y_{s}\}$ is stable set in $G$. Furthermore, we take $S_{2}$ a maximal stable in $G$ such that $x\in S_{2}$. This implies that $S_{2}\setminus x$ is a stable set in $G'$. Then, $|S_{2}|\leq \beta(G')+1\leq \vert S\vert +s=\vert S_{1}\vert$. This is a contradiction, since $S_{2}$ is a maximal stable set and $G$ is unmixed. \QED

\begin{Lemma} \label{7} Let $G$ be an unmixed connected graph without 3-cycles and 5-cycles. If $C$ is a 7-cycle and $H$ is a c-minor of $G$ with $C\subseteq H$ such that $C$ has three non adjacent  vertices of degree 2 in $H$ then $C$ is a c-minor of $G$.
\end{Lemma}

\proof We take a minimal c-minor $H$ of $G$ such that $C\subseteq H$ and $C$ has three non adjacent vertices of degree 2 in $H$. We can suppose that $C=(x,z_{1},w_{1},a,b,w_{2},z_{2})$ with deg$_{H}(x)=$ deg$_{H}(w_{1})$ $=$ deg$_{H}(w_{2})$ $=2$. If $\{z_{1},b\}\in E(H)$ then, $(z_{1},b,w_{2},z_{2},x)$ is a 5-cycle of $G$. Thus, $\{z_{1},b\}\notin E(H)$, similarity $\{z_{2},a\}\notin E(H)$. Furthermore, since $G$ does not have 3-cycles then $\{z_{1},z_{2}\},\{z_{1},a\}$, $\{z_{2},b\}\notin E(H)$. Hence, $C$ is an induced cycle in $H$. On the other hand, if there exists $v\in V(H)$ such that d$(v,C)\geq 2$, then $H'=H\setminus N_{G}[v]$ is a c-minor of $G$ and $C\subseteq H'$. This is a contradiction by the minimality of $H$. Therefore, d$(v,C)\leq 1$ for each $v\in V(H)$.

Now, if deg$_{H}(b)\geq 3$ then there exists $c\in V(H)\setminus V(C)$ such that $\{b,c\}\in E(H)$. If $\{c,z_{2}\}\notin E(G)$ then $N_{H_{1}}(z_{2})$ has two leaves $w_{2}$ and $x$ in $H_{1}=H\setminus N_{H}[w_{1},c]$. This is a contradiction by Lemma \ref{21}. Thus, $\{c,z_{2}\}\in E(H)$. Furthermore, $\{a,c\}$, $\{z_{1},c\}\notin E(H)$ since $(a,b,c)$ and $(z_{1},w_{1},a,b,c)$ are not cycles in $G$. Hence, if deg$_{H}(c)\geq 3$ then there exists $d \in V(H)\setminus V(C)$ such that $\{c,d\}\in E(H)$. Also, $\{d,b\},\{d,z_{2}\},\{d,z_{1}\} \notin E(H)$ since $(c,b,d)$, $(z_{2},b,c)$ and $(z_{1},x,z_{2},c,d)$ are not cycles of $G$. But, d$(d,C)\leq 1$ then, $\{a,d\}\in E(H)$. Consequently, $N_{H_{2}}(z_{1})$ has two leaves $w_{1}$ and $x$ in $H_{2}=H\setminus N_{H}[d,w_{2}]$. This is a contradiction, by Lemma \ref{21} then, deg$_{H}(c)=2$. But, $N_{H_{3}}(z_{2})$ has two leaves $w_{2}$ and $c$ in $H_{3}=H\setminus N_{H}[a]$. A contradiction, therefore, deg$_{H}(b)=2$. Similarity, deg$_{H}(a)=2$.

Now, if deg$_{H}(z_{2})\geq 3$ then  there exists $c'\in V(H)\setminus V(C)$ such that $\{c',z_{2}\}\in E(H)$. If there exists $d'\in V(H)\setminus V(C)$ such that $\{c',d'\}\in E(H)$ then, $\{d',z_{1}\}$ or $\{d',z_{2}\}\in E(G)$, since d$(d',C)\leq 1$. But $(c',d',z_{2})$ and $(x,z_{2},c',d',z_{1})$ are not cycles of $H$. Thus, $N_{H}(c')\subseteq \{z_{1},z_{2}\}$. Consequently, $N_{H_{4}}(z_{2})$ has two leaves $x$ and $c'$ in $H_{4}=H\setminus N_{H}[w_{1}]$. This is a contradiction by Lemma \ref{21}, hence, deg$_{H}(z_{2})=2$. Similarity, deg$_{H}(z_{1})=2$. Furthermore, by the minimality, $H$ is connected. Therefore, $H=C$ and $C$ is a c-minor of $G$. \QED

\section{König and Cohen-Macaulay graphs without 3-cycles and 5-cycles.}

\begin{Proposition} \label{2} Let $G$ be a König graph. If $G'=G\setminus Z_{G}$ then, the following properties are equi\-va-lent:
\begin{description}
 \item{\rm (i)} $G$ is unmixed vertex decomposable.
 \item{\rm (ii)} $\Delta_G$ is pure shellable.
 \item{\rm (iii)} $R/I(G)$ is Cohen-Macaulay.
 \item{\rm (iv)} $G'=\emptyset$ or $G'$ is an unmixed graph with a perfect matching $e_{1},...,e_{g}$ of König type  without squares with two $e_{i}$'s.
 \item{\rm (v)} $G'=\emptyset$ or there exists relabelling of vertices $V(G')=\{x_{1},...,x_{h},y_{1},...,y_{h}\}$ such that $\{x_{1},y_{1}\},...,\{x_{h},y_{h}\}$ is a perfect matching, $X=\{x_{1},...,x_{h}\}$ is a minimal vertex cover of $G'$ and the following conditions holds:
\begin{description}
 \item{\rm (a)} If $a_{i}\in \{x_{i},y_{i}\}$ and $\{a_{i},x_{j}\}, \{y_{j},x_{k}\}\in E(G')$, then $\{a_{i},x_{k}\}\in E(G')$ for distinct $i, j, k$;
 \item{\rm (b)} If $\{x_{i},y_{j}\}\in E(G')$, then $i\leq j$.
\end{description} 
\end{description}
\end{Proposition}

\proof (i)$\Leftrightarrow$(ii)$\Leftrightarrow$(iii) In each case $G$ is unmixed and König then by Lemma \ref{12}, $G$ is totally disconnected or $G'$ is very well-covered. If $G$ is totally disconnected then we obtain the equivalences. Now, if $G'$ is very well covered then by (\cite{MM}, Theorem 1.1) we obtain the equivalences.

(iv)$\Leftrightarrow$(iii) We can assume that $G'\neq \emptyset$ then by Lemma \ref{12}, $G'$ is very well-covered. Hence, by (\cite{half}, Theorem 3.4) since $G=G'\cup Z_{G}$ therefore, $G$ Cohen-Macaulay. 

(iii)$\Rightarrow$ (v) Since $R/I(G)$ is Cohen-Macaulay then $G$ is unmixed. Thus, by Lemma \ref{12}, we can assume that $G'$ is very well-covered and by (\cite{MM}, Lemma 3.1) $G'$ satisfies (v).

(v)$\Rightarrow$(iv) We can assume that $G'\neq \emptyset$. Since, $G'$ has a perfect matching $e_{1},...,e_{h}$ then $\nu(G')=h$. Furthermore, $X$ is a minimal vertex cover, then $\tau(G')=h$ and $e_{1},...,e_{h}$ is a perfect matching of König type. Thus, from (a) and Proposition \ref{22} $G'$ is unmixed. Finally, from (b) there are not square with two $e_{i}$'s. Since, $\{y_{1},...,y_{h}\}$ is a stable set. Therefore, $G'$ satisfies (iv). \QED

\begin{Corollary} \label{cuadrado} Let $G$ be a connected König graph. If $G$ is a Cohen-Macaulay graph then $G$ is an isolated vertex or $G$ has at least a free vertex. \end{Corollary}

\proof By Proposition \ref{2} (v), $G$ is an isolated vertex or $G$ has a perfect matching $e_{1}=\{x_{1},y_{1}\},...,e_{h}=\{x_{h},y_{h}\}$ where $\{x_{1},...,x_{h}\}$ is a minimal vertex cover. Thus, $\{y_{1},...,y_{g}\}$ is a maximal stable set. Furthermore, if $\{x_{i},y_{j}\}$ then $i\leq j$. Therefore, $N_{G}(y_{1})=\{x_{1}\}$, that is, $y_{1}$ is a free vertex. \QED

\begin{Lemma} \label{6} Let $G$ be an connected Cohen-Macaulay graph with a perfect mat\-ching $e_{1}$, ..., $e_{g}$ of König type without square with two $e_{i}$'s and $g\geq 2$. For each $z\in V(G)$ we have that: 

\begin{description}
 \item{\rm (a)} If $\deg_{G}(z)\geq 2$ then there exist $\{z,w_{1}\},\{w_{1},w_{2}\}\in E(G)$ such that $\deg_{G}(w_{2})$ $=1$. Furthermore, $e_{i}=\{w_{1},w_{2}\}$ for some $i\in \{1,...,g\}$.
 \item{\rm (b)} If $\deg_{G}(z)=1$ then there exist $\{z,w_{1}\},\{w_{1},w_{2}\},\{w_{2},w_{3}\}\in E(G)$ such that $\deg_{G}(w_{3})=1$. Furthermore, $e_{i}=\{z,w_{1}\}$ and $e_{j}=\{w_{2},w_{3}\}$ for some $i,j\in \{1,...,g\}$.
\end{description}

\end{Lemma}

\proof Since $e_{1}=\{x_{1},y_{1}\},...,e_{g}=\{x_{g},y_{g}\}$ is a perfect matching of König type we can assume $D=\{x_{1},...,x_{g}\}$ a minimal vertex cover. Thus, $F=\{y_{1},...,y_{g}\}$ is a maximal stable set. We take any vertex $z\in V(G)$. By Proposition \ref{2} we can assume that if $\{x_{i},y_{j}\}\in E(G)$ then $i\leq j$.

(a) First, we suppose that $z=x_{k}$ and there is a vertex $x_{j}$ in $N_{G}(x_{k})$. If $y_{j}$ is a free vertex then, we take $w_{1}=x_{j}$ and $w_{2}=y_{j}$. Consequently, $e_{j}=\{w_{1},w_{2}\}$. Now, we can assume $N_{G}(y_{j})\setminus x_{j}=\{x_{p_1},...,x_{p_r}\}$ with $p_{1}<\cdots < p_{r}< j$. If $y_{p_1}$ is not a free vertex then, there is a vertex $x_{p}$ with $p<p_{1}$ such that $\{x_{p},y_{p_1}\}\in E(G)$. Since $G$ is unmixed and from Proposition \ref{22} we obtain that $\{x_{p},y_{j}\}=(\{x_{p},y_{p_1}\}\setminus y_{p_1})\cup (\{y_{j},x_{p_1}\}\setminus x_{p_1})\in E(G)$. But $p<p_{1}$, this is a contradiction since $p_{1}$ is minimal. Consequently, deg$_{G}(y_{p_1})=1$. Also, since $G$ is unmixed and from Proposition \ref{22} we have that $\{x_{k},x_{p_1}\}=(\{x_{k},x_{j}\}\setminus x_{j})\cup (\{x_{p_1},y_{j}\}\setminus y_{j})\in E(G)$. Therefore, we take $w_{1}=x_{p_1}$ and $w_{2}=y_{p_1}$, furthermore, $e_{p_1}=\{w_{1},w_{2}\}$. Now, we assume that $z=x_{k}$ and $N_{G}(x_{k})\setminus y_{k}=\{y_{j_1},...,y_{j_t}\}$ with $k < j_{1}< \cdots < j_{t}$. We suppose that deg$_{G}(x_{j_t})\geq 2$. If there is a vertex $y_{r}$ such that $\{x_{j_t},y_{r}\}\in E(G)$ then $\{x_{k},y_{r}\}=(\{x_{k},y_{j_t}\}\setminus y_{j_t}) \cup (\{y_{r},x_{j_t}\}\setminus x_{j_t}) \in E(G)$, since $G$ is unmixed and Proposition \ref{22}. This is a contradiction, since $r>j_{t}$ and $j_{t}$ is maximal. Thus, there exist a vertex $x_{p}$ such that $\{x_{j_t},x_{p}\}\in E(G)$. But, since $G$ is unmixed then, $\{x_{k},x_{p}\}=(\{x_{k},y_{j_t}\}\setminus y_{j_t})\cup(\{x_{p},x_{j_t}\}\setminus x_{j_t})\in E(G)$ by Proposition \ref{22}. This is a contradiction, since $N_{G}(x_{k})\setminus y_{k}=\{y_{j_1},...,y_{i_t}\}$. Hence, deg$_{G}(x_{j_t})=1$. Therefore, we take $w_{1}=y_{j_t}$ and $w_{2}=x_{j_t}$, implying, $e_{j_t}=\{w_{1},w_{2}\}$.

Finally, we assume that $z=y_{k}$, since $y_{k}$ is not a free vertex then $N_{G}(y_{k})\setminus x_{k}=\{x_{j_1},...,x_{j_r}\}$ with $j_{1}< \cdots <j_{r}<k$. If $y_{j_1}$ is not a free vertex then there is a vertex $x_{q}$ such that $\{x_{q},y_{j_1}\}\in E(G)$ with $q<j_{1}$. This implies $\{x_{q},y_{k}\}=(\{x_{q},y_{j_1}\}\setminus y_{j_1})\cup (\{x_{j_1},y_{k}\}\setminus x_{j_1})\in E(G)$. But $q<j_{1}$, this is a contradiction. Therefore, deg$_{G}(y_{j})=1$ and we take $w_{1}=x_{j_1}$ and $w_{2}=y_{j_1}$, furthermore, $e_{j_1}=\{w_{1},w_{2}\}$.

(b) Since $e_{1},...,e_{g}$ is a perfect matching then there exists $i\in \{1,...,g\}$ such that $e_{i}=\{z,z'\}$. Furthermore, since $G$ is connected, $g\geq 2$ and $z$ is a free vertex, then $\deg_{G}(z')\geq 2$. Hence, by incise (a)  there exist $w_{1}',w_{2}'\in V(G)$ such that $\{z',w_{1}'\},\{w_{1}',w_{2}'\}\in E(G)$ where deg$_{G}(w_{2}')=1$ and $\{w_{1}',w_{2}'\}=e_{j}$ for some $j\in \{1,...,g\}$.  Therefore, we take $w_{1}=z'$, $w_{2}=w_{1}'$, $w_{3}=w_{2}'$. Consequently $e_{i}=\{z,w_{1}\}$. \QED

\begin{Theorem} \label{P1} Let $G$ be a graph without 3-cycles and 5-cycles. If $G_{1},...,G_{k}$ are the connected components of $G$, then the following condition are equivalent:  
\begin{description}
 \item{\rm (a)} $G$ is unmixed vertex decomposable.
 \item{\rm (b)} $G$ is pure shellable.
 \item{\rm (c)} $G$ is Cohen-Macaulay
 \item{\rm (d)} $G$ is unmixed and if $G_{i}$ is not a isolated vertex then $G_{i}$ has a perfect matching $e_{1},...,e_{g}$ of König type without squares with two $e_{i}'s.$
\end{description}
\end{Theorem}

\proof (a) $\Rightarrow$ (b) $\Rightarrow$ (c) It is known 

(d) $\Rightarrow$ (a) Each component $G_{i}$ is König then $G$ is König. Therefore, from Proposition \ref{2}, $G$ is unmixed vertex decomposable. 

(c) $\Rightarrow$ (d) Since $G$ is Cohen-Macaulay then, $G$ is unmixed. Now, by induction on $\vert V(G)\vert$. We take $x\in V(G)$ such that deg$_{G}(x)$ is minimal. If $r=$ deg$_{G}(x)$ then, we will prove that $r\leq 1$. We suppose $N_{G}(x)=\{z_{1},...,z_{r}\}$ with $r\geq 2$. Furthermore, since $G$ does not contain 3-cycles $N_{G}(x)$ is a stable set. By Remark \ref{27}, $G'=G\setminus N_{G}[x]$ is a Cohen-Macaulay graph. We take $G_{1}',...,G_{s}'$, the connected components of $G'$. By induction hypothesis $G'$ satisfies (d). We can assume that $V(G_{i}')=\{y_{i}\}$ for $i\in \{1,...,s\}$. Since deg$_{G}(x)$ is minimal then, $\{y_{i},z_{j}\}\in E(G)$ for all $i\in \{1,...,s\}$ and $j\in \{1,...,r\}$. If $k=k'$ then, the only maximal stable sets of $G$ are $\{y_{1},...,y_{k},x\}$ and $\{z_{1},...,z_{r}\}$. Thus, $G$ is bipartite, hence $G$ has a free vertex. This is a contradiction. 
Consequently, there is a component $G_{i}'$ with an edge $e=\{w,w'\}$. Since deg$_{G}(x)$ is minimal then, there exist $a,b\in V(G)$ such that $\{a,w\},\{b,w\}\in E(G)$. If $a=b$ then, $(a,w,w')$ is a 3-cycle in $G$. Hence, $a\neq b$. Now, if $a,b\in N_{G}(x)$ then $(x,a,w,w',b)$ is a 5-cycle in $G$. Thus, $|\{w,w',a,b\}\cap V(G_{i}')|\geq 3$. Furthermore, since $G_{i}'$ has a perfect matching then, $\tau(G_{i}')\geq 2$. Furthermore, by Corollary \ref{cuadrado}, $G_{i}'$ has a free vertex $a'$. Then, by Lemma \ref{6} (b), there exist edges $\{a',w_{1}\},\{w_{1},w_{2}\},\{w_{2},b'\}\in E(G_{i}')$ such that deg$_{G_{i}'}(a')=$ deg$_{G_{i}'}(b')=1$. By the minimality of deg$_{G}(x)$ we have that $a$ and $b$ are adjacent with at least $r-1$ neighbour vertices of $x$. If $r\geq 3$ then there exists $z_{j}$ such that $z_{j}\in N_{G}(a)\cap N_{G}(b)$. This implies, $(a,w_{1},w_{2},b,z_{j})$ is a 5-cycle of $G$, this is a contradiction. Consequently, $r=2$. Furthermore, we can assume that $\{a,z_{1}\}, \{b,z_{2}\}\in E(G)$, implying $C=(x,z_{1},a,w_{1},w_{2},b,z_{2})$ is a 7-cycle with deg$_{G}(b)=$ deg$_{G}(a)=$ deg$_{G}(x)=2$. Hence, by Lemma \ref{7}, $G$ has $C$ as a c-minor. This is a contradiction, since a 7-cycle is not Cohen-Macaulay and the c-minors of Cohen-Macaulay graphs are Cohen-Macaulay graphs. Therefore, deg$_{G}(x)=r\leq 1$. 

If $r=0$ then, the result is clear. Now, if $r=1$ then, $N_{G}(x)=\{z\}$. We can assume that $z\in V(G_{1})$. Thus, the connected components of $G\setminus N_{G}[x]$ are $F_{1},G_{2},...,G_{k}$ where $F_{1}=G_{1}\setminus N_{G}[x]$. By induction hypothesis $G_{2},...,G_{k}$ satisfy (d) and there exist $e_{1}=\{x_{1},y_{1}\},...,e_{g}=\{x_{g'},y_{g'}\}$ a perfect matching  of König type of $F_{1}$. We can assume that $D'=\{x_{1},...,x_{g'}\}$ is a minimal vertex cover of $G'$. Consequently, $e_{1},...,e_{g'},e_{g'+1}=\{x,z\}$ is a perfect matching of $G_{1}$ . Furthermore, $D=\{x_{1},...,x_{g'},z\}$ is a vertex cover of $G$, then $g'+1\leq \nu(G_{1}) \leq \tau(G_{1}) \leq |D|$. Hence, $\nu(G_{1})=\tau (G_{1})=g'+1$ and $e_{1},...,e_{g'+1}$ is a perfect matching of König type of $G_{1}$. Furthermore, there are not square with two $e_i$'s in $e_{1},...,e_{g'}$ and $x$ is a free vertex. Therefore, $G_{1}$ does not contain square with two $e_{i}$'s for $i=1,...,g'+1$. \QED

\begin{Corollary} Let $G$ be a graph without 3-cycles and 5-cycles. If $G$ is Cohen-Macaulay then $G$ has at least an extendable vertex $x$ adjacent to a free vertex.
\end{Corollary}
 
\proof From Theorem \ref{P1}, $G$ is König. Thus, by Corollary \ref{cuadrado} there exists a free vertex $x$. If $N_{G}(x)=\{y\}$, therefore, by Remark \ref{25}, $y$ is an extendable vertex. \QED

\begin{Definition} Let $G$ be a graph and let $S\subseteq V(G)$. For use below consider the graph $G\cup W_{G}(S)$ obtained from $G$ by adding new vertices $\{y_{i}\ |\ x_{i}\in S\}$ and new edges $\{\{x_{i},y_{i}\}\ |\ x_{i}\in S\}$. The edges $\{x_{i},y_{i}\}$ are called \textit{whiskers} (pendants). $G$ is called \textit{whisker graph} if there exists an induced subgraph $H$ of $G$ such that $V(H)=\{x_{1},...,x_{s}\}$, $V(G)=V(H)\cup \{y_{1},...,y_{s}\}$ and $E(G)=E(H)\cup W(H)$ where $W(H)=\{\{x_{1},y_{1}\},...,\{x_{s},y_{s}\}\}$. Its pendant edges form a perfect matching. \end{Definition}

\begin{Corollary} If $G$ is a connected graph of girth 6 or more, then the following condition are equivalent:

\begin{description}
 \item{\rm (i)} $G$ is unmixed vertex decomposable.
 \item{\rm (ii)} $\Delta_G$ is pure shellable.
 \item{\rm (iii)} $R/I(G)$ is Cohen-Macaulay.
 \item{\rm (iv)} $G$ is unmixed and König.
 \item{\rm (v)}  $G$ is an isolated vertex or $G$ very well-covered.
 \item{\rm (vi)}  $G$ is unmixed with $G\neq C_{7}$.
 \item{\rm (vii)}  $G$ is an isolated vertex or $G$ unmixed whiskers.
\end{description}

\end{Corollary}

\proof (i) $\Rightarrow$ (ii) $\Rightarrow$ (iii) It is known. (iii) $\Rightarrow$ (iv) $G$ is unmixed and from Theorem \ref{P1}, $G$ is König. (iv) $\Rightarrow$ (v) From Lemma \ref{12}. (v) $\Rightarrow$ (vi) It is clear, since $C_{7}$ is not very well-covered. 

(vi) $\Rightarrow$ (vii) By (\cite{Finbow}, Corollary 5), if $G$ is not an isolated vertex then, the pendant edges $\{x_{1},y_{1}\},...,\{x_{g},y_{g}\}$ of $G$ form a perfect matching since $\{x_{i},y_{i}\}$ is a pendant edge. We can assume that $\deg_{G}(y_{i})=1$ for all $i= 1,...,g$. Hence, if $G'=G[x_{1},...,x_{n}]$ and $G=G'\cup W(V_{G'})$ with $W(V(G'))=\{y_{1},...,y_{g}\}$, therefore $G$ is the whisker graph. 

(vii) $\Rightarrow$ (i) If $G$ is an isolated vertex, it is clear. Now, if $G$ is a whisker graph, then there exists a perfect matching $e_{1}=\{x_{1},y_{1}\},...,e_{g}=\{x_{g},y_{g}\}$ such that deg$_{G}(y_{i})=1$ for $i=1,...,g$. Thus, $D=\{x_{1},...,x_{g}\}$ a minimal vertex cover and $\tau(G)=g$. Hence, $e_{1},...,e_{g}$ is a perfect matching of König type without square with two $e_{i}$'s. Therefore, by Theorem \ref{P1}, $G$ is unmixed vertex decomposable. \QED

\section{Unicyclic graphs.}

\begin{Definition} A \textit{unicyclic} graph is a connected graph with exactly one cycle. \end{Definition}

\begin{Theorem} \label{9} Let $G$ be a unicyclic graph whose cycle is $C$ then, the following conditions are equivalent: 
\begin{description}
 \item{\rm (1)} $G$ is vertex decomposable.
 \item{\rm (2)} $G$ is shellable
 \item{\rm (3)} $G$ is sequentially Cohen-Macaulay.
 \item{\rm (4)} $\vert V(C) \vert\in \{3,5\}$ or there exists $x\in D_{1}(C)$ such that \rm deg$_{G}(x)=1$
\end{description} 
\end{Theorem} 

\proof (1) $\Rightarrow$ (2) $\Rightarrow$ (3)  It is known. 

(3) $\Rightarrow$ (4) If $\vert V(C) \vert\notin \{3,5\}$ then, the cycle $C$ is not sequentially Cohen-Macaulay. Furthermore, $C$ is the only 2-connected block in $G$ then, by Theorem \ref{5} there exists $x\in D_{1}(C)$ such that deg$_{G}(x)=1$.  

(4) $\Rightarrow$ (1) If $\vert V(C) \vert =\{3,5\}$ then, $G$ is vertex decomposable by (\cite{RW}, Theorem 1) 

Now, if $\vert V(C)\vert \neq 3,5$, then $C$ is not vertex decomposable. Hence, by Corollary \ref{3} there exists $x\in D_{1}(C)$ such that deg$_{G}(x)=1$. If $N_{G}(x)=\{y\}$ then $y$ is a shedding vertex. Since $y\in V(C)$, then $G\setminus y$ and $G\setminus N_{G}[y]$ are forests, consequently they are vertex decomposable. Therefore, $G$ is vertex decomposable. \QED

\begin{Definition} A clique $S$ of a graph $G$ containing at least one simplicial vertex of $G$ is called a \textit{simplex} of $G$. \end{Definition}

\begin{Definition} A 4-cycle is called  \textit{basic} if it contains two adjacent vertices of degree two, and the remaining two vertices belong to a simplex or a basic 5-cycle of $G$. \end{Definition}

\begin{Remark} In the case unicyclic graph we have that a 4-cycle is basic if it contains two adjacent vertices of degree two and the remaining two vertices is jointed to free vertex.
\end{Remark}

\begin{Definition} A graph $G$ is in the family $\mathcal{SQC}$, if $V(G)$ can be partitioned into three disjoint subsets $S_{G}$, $Q_{G}$ and $C_{G}$: the subset $S_{G}$ contains all vertices of the simplexes of $G$, and the simplexes of $G$ are vertex disjoint; the subset $C_{G}$ consists of the vertices of the basic $5$-cycles and the basic 5-cycles form a partition of $C_{G}$; the remaining set $Q_{G}$ contains all vertices of degree two of the basic 4-cycles.  \end{Definition}

\begin{Remark} \label{31} If $G\in \mathcal{SQC}$ then $G$ is vertex decomposable by {\rm(\cite{ON}, Theorem 2.3)} and $G$ is unmixed by {\rm(\cite{cactus}, Theorem 3.1)}. \end{Remark}

\begin{Lemma} \label{33} If $x$ is a free vertex of $G$ such that $G'=G\setminus N_{G}[x]$ is an unmixed whiskers. If there are not cycle $C$ of $G$ such that $z\in V(C)$ where $\{x,z\}\in E(G)$ then $G$ is a whisker.
\end{Lemma}

\proof Let $G_{1},...,G_{s}$ be the connected components of $G$ then there exist induced subgraphs $H_{1},...,H_{s}$ such that $V(H_{i})=\{x_{1}^{i},...,x_{r_{i}}^{i}\}$, $V(G_{i})=V(H_{i})\cup \{y_{1}^{i},...,y_{r_{i}}^{i}\}$ and $E(G_{i})=E(H_{i})\cup W(H_{i})$ where $W(H_{i})=\{\{x_{1}^{i},y_{1}^{i}\},...,\{x_{r_{i}}^{i},y_{r_{i}}^{i}\}\}$. If $\{x,z\}$ is a connected component of $G$ then $G=\{x,z\}\cup G_{1}\cup \cdots \cup G_{s}$ and each component is a whisker graph, therefore, $G$ is a whisker graph. Now, if there are $w_{1}^{i},w_{2}^{i}\in E(G_{i})$ such that $\{w_{1}^{i},z\},\{w_{2}^{i},z\} \in E(G)$ then there exists a path $\{w_{1}^{i},v_{1},...,v_{m},w_{2}^{i}\}$ since $G_{i}$ is connected and $v_{k}\in V(G_{i})$ for $k\in \{1,...,m\}$, therefore, $z\in V(C)=(z,w_{1}^{i},v_{1},...,v_{m},w_{2}^{i})$, this is a contradiction by hypothesis. Thus $\vert N_{G}(z)\cap V(G_{i})\vert\leq 1$ for $i\in \{1,...,s\}$. If $\vert V(G_{i})\vert =2$ then $\deg_{G_{i}}(x_{1}^{i})=\deg_{G_{i}}(y_{1}^{i})=1$, therefore, if $\{w,z\}\in E(G)$ with $w\in \{x_{1}^{i},y_{1}^{i}\}$ then $G_{i}\cup \{x,z\}$ is a whisker graph. We can suppose $\vert V(G_{i})\vert\geq 3$, then if $\{z,y_{j}\}\in E(G)$ for some $j\in \{1,...,r_{i}\}$ then there is a vertex $w$ such that $\{w,x_{r}\}\in V(G_{i})$ since $\vert V(G_{i}) \vert \geq 3$ then $N_{G'}(z)$ has two leaves $y_{j}$ and $x$ in $G'=G\setminus N_{G}[w]$, this is contradiction, therefore, $\{z,y_{j}\}\notin E(G)$ for all $j\in \{1,...,r_{i}\}$. Thus, $\{x_{j}^{i},z\}\in E(G)$ for some $j\in \{1,...,r_{i}\}$, we can take $V(H)=\bigcup_{i=1}^{s}V(H_{i})\cup \{z\}$, $V(G)=\bigcup _{i=1}^{s}V(G_{i})\cup \{z,x\}$ and $E(G)=E(H)\cup W(H)$ where $W(H)=\bigcup_{i=1}^{s}W(H_{i})\cup \{x,z\}$. Therefore, $G$ is a whisker graph. \QED

\begin{Remark} \label{32} In {\rm(\cite{Villa}, Theorem 7.3.17)} is proved that if $G$ is a tree then $G$ is a Cohen-Macaulay graph if and only if $G$ is unmixed if and only if $G$ is a whisker graph. in the following result we characterize the unicyclic graph unmixed and Cohen-Macaulay.
\end{Remark}

\begin{Theorem} \label{23} Let $G$ be a connected unicyclic graph whose cycle is $C$ then, $G$ is well-covered if and only if $G$ satisfies one of the following conditions: 
\begin{description}
 \item{\rm (a)} $G\in \{C_{3},C_{4},C_{5}C_{7}\}$. 
 \item{\rm (b)} $G$ is a whisker graph.
 \item{\rm (c)} $C$ is a simplex 3-cycle or a basic 5-cycle and $G\setminus V(C)$ is a whisker graph. 
 \item{\rm (d)} $C$ is a basic 4-cycle with two adjacent vertices $a,b$ of degree 2 in $G$ such that $G\setminus \{a,b\}$ is a whisker graph. 
\end{description}
\end{Theorem}

\proof $\Rightarrow$) By induction on $l=\vert V(G)\setminus V(C)\vert$. If $l=0$ then, $G=C$ and $G$ satisfies (a) since the unmixed cycles are $C_{3},C_{4},C_{5}$ and $C_{7}$. Now, if $l\geq 1$ then, we can suppose that $C=(y_{1},...,y_{k})$ with $k \geq 3$. Hence, $C\subsetneq G$. Since $C$ is the only cycle in $G$ then, there is a free vertex $x$ such that $N_{G}(x)=\{z\}$. We can take $x$ such that d$(x,C)=max\{d(a,C) \vert\ a\in V(G),\ \deg_{G}(a)=1\}$. Furthermore, we can suppose that $N_{G}(x)=\{z\}$.  Since $G$ is unmixed then, $G'=G\setminus N_{G}[x]$ is unmixed by Remark \ref{27}. 

If d$(x,C)\geq 2$ then $C\subseteq G'$ and $G'$ satisfies (a), (b), (c) or (d) by induction hypothesis. If $G'$ satisfies (a) then, $G'=C_{r}$ with $r\in \{3,4,5,7\}$. We can assume that $C_{r}=(y_{1},...,y_{r})$ and $\{z,y_{1}\}\in E(G)$. If $r=4$ or $7$ then $N_{G_{2}}(z)$ has two leaves $x$ and $y_{1}$ in $G_{2}=G\setminus N_{G}[y_{3},y_{r-1}]$. This is a contradiction by Lemma \ref{21}, if $r=3$ or $5$. Therefore, $G_{1}'=G'\setminus N_{G'}[y_{2}]=V(C)$ is the edge $\{x,z\}$. Then $G_{1}'\cup \{x,z\}$ is a whisker graph and $C_{3}$ or $C_{5}$ is a simplex 3-cycle or a basic 5-cycle in $G$. Hence, $G$ satisfies (c). If $G'$ satisfies (b) then $G$ is a whisker graph by Lemma \ref{33}. Hence, $G$ satisfies (b). If $G'$ satisfies (c) then $C$ is a simplex 3-cycle or a basic 5-cycle and $G_{1}'$ is a whisker forest graph. $C_{3}$ is simplex in $G$ since if $\{z,y_{1}\}\in E(G)$ and there are not simplicial vertices in $C_{3}$ then there are $z_{2},z_{3}\in V(G)$ such that $\{y_{2},z_{2}\},\{y_{3},z_{3}\}\in E(G)$. Therefore, $N_{G_{3}}(z)$ has two leaves in $G_{3}=G\setminus N_{G}[z_{2},z_{3}]$. This is a contradiction. If $C_{5}$ is not basic in $G$ then there are two adjacent vertices in $C_{5}$ such that $\deg_{G}(y_{1})\geq 3$ and $\deg_{G}(y_{2})\geq 3$ then there is $z_{2}$ such that $\{y_{2},z_{2}\}\in E(G)$. Therefore, $N_{G_{4}}(z)$ has two leaves $y_{1}$ and $x$ in $G_{4}=G\setminus N_{G}[z_{2},y_{4}]$. This is a contradiction. Then, $C_{5}$ is a basic 5-cycle in $G$. Furthermore, $G_{1}=G_{1}'\cup \{x,z\}$ is a whisker graph by Lemma \ref{33}. Thus, $G$ satisfies (c). If $G'$ satisfies (d). Hence, $C$ is a basic 4-cycle $(a,b,c,d)$ with $\deg_{G'}(a)=\deg_{G'}(b)=2$ such that $G'\setminus \{a,b\}$ is a whisker graph. If $\{a,z\}\in E(G)$ then $N_{G_{5}}(z)$ has two leaves $x$ and $a$ in $G_{5}=G\setminus N_{G}[c]$, this is a contradiction by Lemma \ref{21}. Then $\{a,z\}\notin E(G)$. Similarity, $\{b,z\}\notin E(G)$. Therefore, deg$_{G}(a)=$deg$_{G}(b)=2$ and $G\setminus \{a,b\}=(G'\setminus \{a,b\})\cup \{x,z\}$ is a whisker graph by Lemma \ref{33}. Thus, $G$ satisfies (d).  

If d$(x,C)=1$ then $V(G)=V(C)\cup D_{1}(C)$ and $G'$ is a forest, therefore, $G'$ is a whisker graph by Remark \ref{32}, also, there is a vertex $y_{i}\in V(C)$ such that $z=y_{i}$. We can assume $i=1$. Since $G$ is unicyclic, if $w\in D_{1}(C)$ then $w$ is a free vertex. If each $y_{i}\in V(C)$ is adjacent to one free vertex then $G$ is a whisker graph. Then $G$ satisfies (b). If there is a vertex $y_{j}\in V(C)$ such that deg$_{G}(y_{j})=2$, without loss of generality, we can assume $i=2$ then, $\vert V(C)\vert\leq 5$ by (\cite{Finbow}, Corollary 5). If $r=3$ then $C$ is a simplex 3-cycle. Furthermore, $G_{1}=G\setminus V(C)=G\setminus N_{G}[y_{2}]$ is not a whisker graph since the vertex $\{x\}$ is isolated in $G_{1}$. Therefore, $\deg_{G}(y_{i})=3$ for all $i\in \{1,2,3\}$. Thus, $G$ satisfies (b). Now, if $r=5$, then $N_{G_{7}}(y_{1})$ has two leaves $x$ and $y_{2}$ in $G_{7}=G\setminus N_{G}[y_{4}]$. This is not possible. Therefore, each $y_{i}\in V(C)$ is adjacent to a free vertex. Hence, $G$ is a whisker graph. Thus, $G$ satisfies (b). If $r=4$, if deg$_{G}(y_{4})=2$ then $N_{G'}(y_{3})$ has two leaves $y_{2}$ and $y_{4}$ in $G'$. Consequently,  deg$_{G}(y_{4})\geq 3$ and there is a free vertex $z_{4}$ such that $\{y_{4},z_{4}\}\in E(G)$. Now, if deg$_{G}(y_{3})\geq 3$, then there is a free vertex $z_{3}$ such that $\{y_{3},z_{3}\}\in E(G)$. Thus, $N_{G_{8}}(y_{3})$ has two leaves $y_{2}$ and $z_{3}$ in $G_{8}=G\setminus N_{G}[x,z_{4}]$. Hence, $\deg_{G}(y_{3})=2$ and $C$ is a basic 4-cycle. Furthermore, $G\setminus \{y_{2},y_{3}\}$ is a whisker graph. Therefore, $G$ satisfies (d). 
 
$\Leftarrow$) The families (a),(b),(c) and (d) are in $\mathcal{SQC} \cup \{C_{4},C_{7}\}$, then by Remark \ref{31} they are unmixed. \QED

\begin{Corollary} \label{13} Let $G$ be an unicyclic graph then, the following condition are equivalent:
\begin{description}
 \item{\rm (1)} $G$ is unmixed vertex decomposable.
 \item{\rm (2)} $G$ is pure shellable
 \item{\rm (3)} $G$ is Cohen-Macaulay.
 \item{\rm (4)} $G$ is unmixed and $G\neq C_{4}, C_{7}$.
\end{description}
\end{Corollary} 

\proof (1)$\Leftrightarrow$(2)$\Leftrightarrow$(3) By Theorem \ref{9}.

(3) $\Rightarrow$ (4) Since $G$ is Cohen-Macaulay then $G$ is unmixed. Furthermore, $C_{4}$ and $C_{7}$ are not Cohen-Macaulay. 

(4) $\Rightarrow$ (1) Let $C$ be the cycle of $G$. By Theorem \ref{23}, $G\in \mathcal{SQC}$ since $G\neq C_{4},C_{7}$. Therefore $G$ is vertex decomposable by Remark \ref{31}. \QED

\begin{Corollary} If $G$ is an unmixed unicyclic then $G$ is vertex decomposable if and only if $G\neq C_{4},C_{7}$.
\end{Corollary}

\proof By Corollary \ref{13}. \QED

\bibliographystyle{plain}

\hspace{4.5mm} Departamento de Matem\'aticas, Centro de Investigaci\'on y de Estudios Avanzados del IPN, Apartado Postal 14-740, 07000 Mexico City, D.F.  \vspace{-2mm}

\textit{E-mail address}: \texttt{idcastrillon@math.cinvestav.mx}

\hspace{4.5mm} Instituto de Matematicas, Universidad de Antioquia, Apartado A\'ereo 1226, Medell\'in, Colombia.
\vspace{-2mm}

\textit{E-mail address}: \texttt{roberto.cruz@udea.edu.co}

\hspace{4.5mm} Departamento de Matem\'aticas, Centro de Investigaci\'on y de Estudios Avanzados del IPN, Apartado Postal 14-740, 07000 Mexico City, D.F.
\vspace{-2mm}

\textit{E-mail address}: \texttt{ereyes@math.cinvestav.mx}

\end{document}